\documentclass[12pt]{amsart}
\usepackage{amssymb,psfig}

\renewcommand{\L}{\mathcal{L}}

\newcommand{\R}{\mbox{$\mathbb{R}$}}

\newcommand{\eps}{\varepsilon}

\newcommand{\M}{\mathcal{M}}

\newcommand{\Z}{\mbox{$\mathbb{Z}$}}

\newcommand{\la}{\lambda}
\newcommand{\al}{\alpha}

\renewcommand{\th}{\theta}
\newcommand{\Om}{\Omega}

\newcommand{\BN}{\mathbb{N}}
\newcommand{\F}{\mathcal{F}}
\newcommand{\ov}{\overline}
\newcommand{\wh}{\widehat}

\newcommand{\om}{\omega}

\newtheorem{lemma}{Lemma}[section]

\newtheorem{prop}[lemma]{Proposition}
\newtheorem{thm}[lemma]{Theorem}
\newtheorem{cor}[lemma]{Corollary}
\theoremstyle{definition}

\newtheorem{exam}[lemma]{Example}

\theoremstyle{remark}
\newtheorem{rmk}[lemma]{Remark}

\begin{document}

\title[On the fine structure of stationary measures]{On the fine structure
of stationary measures in systems which contract-on-average}

\author {Matthew Nicol}
\address {Mathematics Department, University of Surrey, Guildford,
Surrey GU2 5XH, United Kingdom. E-mail: M.Nicol@surrey.ac.uk}
\author {Nikita Sidorov}
\address {Department of Mathematics, UMIST, P.O. Box 88, Manchester M60 1QD,
United Kingdom. E-mail: Nikita.A.Sidorov@umist.ac.uk}
\author {David Broomhead}
\address {Department of Mathematics, UMIST, P.O. Box 88, Manchester M60 1QD,
United Kingdom. E-mail: David.Broomhead@umist.ac.uk}

\thanks{This work was made possible
by  EPSRC grant number GR/L98923.
MJN supported in part by the LMS and the Nuffield Foundation.
We would like to thank Vadim Kaimano\-vich, Anders \"Oberg, Yuval Peres,
Mark Pollicott, Vladimir Re\-meslen\-nikov and Boris Solomyak
for helpful conversations, communications and suggestions.}
\date{\today}

\begin{abstract}
Suppose $\{ f_1,\dots,f_m\}$ is a set of Lipschitz maps
of $\R^d$. We form the iterated function system (IFS)
by independently choosing the maps so that the map  $f_i$ is chosen
with  probability $p_i$  ($\sum_{i=1}^m p_i=1$).
We assume that the IFS contracts on average.
We give an upper bound for the Hausdorff dimension of
the invariant measure induced on $\R^d$ and as a corollary
show that the measure will be  singular if the modulus of the entropy
$\sum_i p_i \log p_i$ is less  than $d$ times the modulus of the
Lyapunov exponent of the system. Using a version of Shannon's Theorem
for random walks on semigroups we improve this estimate and show
that it is actually attainable for certain cases of affine mappings of $\R$.
\end{abstract}

\maketitle

\section{Introduction}

Suppose $\{f_1,\dots,f_m\}$ is a set of Lipschitz maps
of $\R^d$. We may form an iterated function
system (IFS) by independently choosing the maps
so that the map  $f_i$ is chosen with  probability $p_i$
($\sum_{i=1}^m p_i=1$). We denote the the probability vector by
$\ov p:=(p_1,\dots,p_m)$, and the IFS itself will be denoted by $\Phi$.

More precisely, let $\Omega=\prod_0^\infty\{1,\dots,m\}$
and equip $\Omega$ with
the  product probability measure $\nu$ induced in the standard way
on cylinder sets by the probability vector $\ov p$.
Let  $x_0\in\R^d$. For any $\omega\in\Omega$ and any $n\in\mathbb{N}$
we define the point
$$
x_n(\omega):=f_{\om_0}\dots f_{\om_{n-1}}(x_0).
$$
If $\lim_{n\rightarrow \infty}x_n (\omega)$ exists, then we define
\begin{equation}
\phi(\om)=\lim_{n\to\infty} x_n(\omega).
\label{phi}
\end{equation}

We are going to formulate the hypothesis on $\Phi$ such that the
mapping~$\phi:\Omega \rightarrow \R^d$ will be  defined $\nu$-a.e. Define
\[
h(\ov p):=\sum_{i=1}^m p_i \log p_i.
\]
Note that $-h(\ov p)$ is the measure-theoretic  entropy of the Bernoulli
shift $\sigma:\Omega\to\Omega$ with the probabilities $(p_1,\dots,p_m)$.

For any Lipschitz map $g$  of $\R^d$ we let $\|g\|$ denote the Lipschitz
constant of $g$.
We  assume a  {\it contraction on average}
(sometimes called  {\it  logarithmic average contractivity})
condition to hold: for $\nu$-a.e. $\omega \in \Omega$,
\begin{equation}
\label{contractavcompl}
\lim_{n\to\infty}\frac1n\log^+\|f_{\om_0}\dots f_{\om_{n-1}}\|=\chi(\Phi)<0.
\end{equation}
We will call $\chi(\Phi)$ the {\it Lyapunov exponent} of the system.

Note that the condition~(\ref{contractavcompl}) is implied by the
condition
\[
\sum_{i=1}^m p_i \log \|f_i\|<0.
\]

Contraction on average implies that $\phi$ is a well
defined $\nu$-measurable function defined
by the formula~(\ref{phi}), it is independent of the choice of
initial point $x_0$ and that there exists an
invariant attracting set $A\subset \Omega\times\R^d$
which is the graph of $\phi$. This result is
standard and can be found for instance in \cite[Theorem~3]{El2},
~\cite[Theorem~1.4]{Stark1},~\cite[Theorem~5]{Campbell},
\cite[Theorem~4]{Arnold&Crauel} and
~\cite[Proposition~2.3]{BHN}.

It is well known (see P.~Diaconis and D.~Freedman \cite{DiFr})  that  for any such IFS
there exists a unique {\em stationary measure}
$\mu$ on $\R^d$ independent of the choice of initial point, i.e,
such that $L^*\mu=\mu$, where $L$ is the Perron-Frobenius operator for
the IFS~$\Phi$:
$$
L\psi(x)=\sum_{i=1}^m p_i\psi(f_ix).
$$
In fact the measurable mapping $\phi$ induces
$\mu$ on Borel sets of $\R^d$ by $\mu (B)=\nu\circ\phi^{-1} (B)$.
Sometimes we will also call $\mu$ the {\em invariant measure}.

By results of L. Dubins and D. Freedman \cite{Dubins-Freedman}
(see also M.~Barnsley and J.~Elton \cite[Proposition 1]{mBb88}) on Markov operators,
 $\mu$ must be of {\it pure type},
i.e., either absolutely continuous or purely singular with
respect to Lebesgue measure on $\R^d$. There are also results due
to M.~Barnsley and J.~Elton \cite[Theorem 3]{mBb88} about the structure
of the support of $\mu$ when $d=1$ and the maps $\{ f_i\}$
are affine (see Section~\ref{examples}).

An important classical example of an IFS is the one-parameter family
\begin{align*}
f_0(x)  &  =\la^{-1}x,\\
f_1(x)  &  =\la^{-1}x+1-\la^{-1}
\end{align*}
with $p_1\in(0,1)$. It has been extensively studied since the
1930's. In  recent work by B.~Solomyak \cite{Solomyak} it was shown that
if $p_1=p_2=\frac12$, then a.e. $\lambda\in(1,2)$ induces an absolutely
continuous measure $\mu$ on the interval $[0,1]$. A similar
result was later obtained by the same authors for $p_1\in[1/3,2/3]$
(see Section~\ref{examples}).
However the problem of whether the invariant measure
(usually called the {\em Bernoulli convolutions} or the
{\em Erd\"os measure})
for this system is absolutely continuous or singular for a
{\bf given} value of $\la$ (known as the Erd\"{o}s Problem),
is very hard and only few concrete results are known
(see \cite{PeScSo} for a nice review and collection of
references).

The purpose of this paper is to investigate conditions
on IFS which contract on average  under which their
invariant measure is known to be singular or absolutely continuous.
The structure of the paper is as follows:
in Section~\ref{sec-sing} we present an upper bound for the Hausdorff
dimension of the invariant measure $\mu$ and describe sufficient conditions
for $\mu$ to be singular in terms of $\chi(\Phi),\ h(\ov p)$
and the expansion
rate of the semigroup generated by $\{f_i\}$. In Section~\ref{examples}
we present several examples showing how to apply the main theorem.
Thus, we have reduced the problem of estimating $\dim_H(\mu)$ to certain
combinatorial and algebraic issues concerning the semigroup in
question.

We would like to emphasize that although our results apply to a general
IFS which contracts-on-average, the most interesting case for us will be
the systems in which {\bf not all} of $f_i$ are uniformly contracting, i.e.,
such that the support of $\mu$ is unbounded. One of the reasons
for doing so is that there are some indications that $\mbox{supp}(\mu)$
in this case will be ``less fractal" than for uniformly contracting systems
(see examples below).

\section{sufficient conditions for singularity of the invariant
measure}\label{sec-sing}

Let $G^+$ denote the semigroup generated by the maps $\{f_1,\dots,f_m\}$.
Its elements are all compositions $f_{\om_0}\circ\dots\circ f_{\om_{n-1}}$
for any $n\in\BN$ and $\om_k\in\{1,\dots,m\}$. It is clear that $G^+$ can be
either the free semigroup $\F_m^+$ (if all such compositions are different)
or a proper subsemigroup of $\F_m^+$. Both possibilities can occur
(see Section~\ref{examples}).

Let $D_n$ denote the set of all words of length $n$ in $G^+$.
In other words, $D_n$ is the set of equivalence classes in
$\prod_0^{n-1}\{1,\dots,m\}$, namely:
$$
(\om^*_0,\dots,\om^*_{n-1})\sim(\om'_0,\dots,\om'_{n-1}) \text{ if }
f_{\om^*_0}\circ\dots\circ f_{\om^*_{n-1}}=f_{\om'_0}\circ\dots\circ
f_{\om'_{n-1}}.
$$

From general considerations  the growth of $G^+$ is
exponential, i.e., there exists $\th\in[1,m]$ such that
\begin{equation}
\th=\lim_{n\to+\infty}\sqrt[n]{\#D_n}.
\label{thetaa}
\end{equation}
Indeed, let $d_n=\#D_n$; then $d_{n+k}\le d_nd_k$, because
$(\om_0,\dots,\om_{n-1})\sim(\om'_0,\dots,\om'_{n-1})$ and
$(\om_n,\dots,\om_{n+k-1})\sim(\om'_n,\dots,\om'_{n+k-1})$ imply
$(\om_0,\dots,\om_{n+k-1})\sim(\om'_0,\dots,\om'_{n+k-1})$. Hence, with a
little work,  (\ref{thetaa}) follows. Obviously, $\th\ge1$ and if
$G^+$ is abelian, then $\th=1$.

Let $H_\mu$ denote the entropy of the random walk on the semigroup
$G^+$ with probabilities $\{p_1,\dots,p_m\}$. It is defined as follows: let
$\mu_n$ be the $n$'th convolution of $(p_1,\dots,p_m)$ on $D_n$, i.e.,
\[
\mu_n([\om^*_0,\dots,\om^*_{n-1}])=\sum_{(\om'_0,\dots,\om'_{n-1})\sim
(\om^*_0,\dots,\om^*_{n-1})}\nu(\om_0=\om'_0,\dots,\om_{n-1}=\om'_{n-1}),
\]
where $[\cdot]$ denotes the equivalence class.

We define
$$
H_n:=-\sum_{[y]\in D_n}\mu_n([y])\log \mu_n([y]),
$$
and finally, for $\th>1$,
\begin{equation}
H_\mu:=\lim_{n\to\infty}\frac{H_n}{n\log\th}
\label{Hmu}
\end{equation}
(it is a standard argument that such a limit exists and equals
the infimum of the corresponding sequence). For $\th=1$ we
set $H_\mu:=0$; it is natural, because $H_n\le\log\#D_n$, whence
$\lim_n H_n/n=0$ in this case. By the definition of $H_\mu$ we have
\begin{equation}
0\le H_\mu\le-\frac{h(\ov p)}{\log m}\le1.
\label{entropy}
\end{equation}
We will need a version of Shannon's Theorem for random walks. In
the case of discrete groups it was proved independently by
Y.~Derriennic \cite{Der} and V.~Kaimano\-vich and A.~Vershik \cite{KaiVer}.
We will adapt the proof from \cite{Der} to our ``semigroup" context
(see also \cite[Theorem~1.6.4]{Kai}).

\begin{lemma}\label{shannon} Let $\om\in\Omega$ and
\[
\mathcal{E}_n(\om)=\{\om'\in\Omega\mid (\om_0,\dots,\om_{n-1})
\sim(\om'_0,\dots,\om'_{n-1})\}.
\]
Then for $\nu$-a.e. $\om$,
\[
\lim_{n\to+\infty}\frac{\log\nu\mathcal{E}_n(\om)}{n}
=-H_\mu\log\th.
\]
\end{lemma}

\noindent {\bf Proof:} Let $[\om]_n$ denote the set
of all words of length $n$ equivalent to $(\om_0,\dots,\om_{n-1})$.
We will identify $[\om]_n$ with
$\mathcal{E}_n(\om)$. Hence $\nu(\mathcal{E}_n(\om))=\mu_n([\om]_n)$.

Let $f_n(\om):=-\log\mu_n([\om]_n)$.
By the same reason as in the proof of formula~(\ref{thetaa}),
\[
\mu_{n+k}([\om]_{n+k})\ge \mu_n([\om]_n)\cdot
\mu_k([\om_n,\dots,\om_{n+k-1}])
\]
for any $\om\in\Omega,n\ge1,k\ge1$. Hence
\[
f_{n+k}(\om)\le f_n(\om)+f_k(\sigma^n\om),\quad n,k\ge1,
\]
and by Kingman's Subadditive Ergodic Theorem, there exists the limit
$\ov f(\om)=\lim_n \frac1n f_n(\om)$ for $\nu$-a.e. $\om$ and
\[
\frac1n\int_{D_n} f_n\, d\mu_n\to \ov f,\quad n\to+\infty
\]
as well. It suffices to note that
$\int_{D_n} f_n\, d\mu_n=-\sum_{[y]\in D_n} \mu_n([y])\log\mu_n([y])=
-H_n$ and apply (\ref{Hmu}).\qed

Now we are ready to formulate the main result of this paper.

\begin{thm}\label{thm-2}
Suppose $\Phi$ is an IFS on $\R^d$ which contracts on average. Then
\begin{equation}
\dim_H (\mu) \le -\frac{H_\mu\log\th}{\chi(\Phi)}.
\label{dimension}
\end{equation}
\end{thm}
This has two immediate corollaries:

\begin{cor}\label{nonfree} If
$$
\frac{h(\ov p)\log\th}{\chi(\Phi)}< d\log m,
$$
then $\mu$ is singular, and
$$
\dim_H(\mu)\le\frac{h(\ov p)\log\th}{\chi(\Phi)\log m}<d.
$$
In particular, if $p_1=\dots=p_m=\frac1m$, then
\begin{equation}
\dim_H(\mu)\le\frac{\log\th}{|\chi(\Phi)|}.
\label{dim-th}
\end{equation}
\end{cor}
\noindent {\bf Proof:} follows from (\ref{entropy}) and the fact that
if $\dim_H(\mu)<d$, then $\mu$ is singular.\qed

\begin{cor}\label{cor-1}
The measure $\mu$  is
singular for any $\Phi$ such that
\begin{equation}
d|\chi(\Phi)|>|h(\ov p)|.
\label{chi>h}
\end{equation}
Besides, if (\ref{chi>h}) is satisfied, then
$$
\dim_H(\mu)\le\frac{h(\ov p)}{\chi(\Phi)}<d.
$$
\end{cor}

\begin{rmk} As far as we know, there have been no analogs
of Theorem~\ref{thm-2} in such a general framework. However,
F.~Przytycki and M.~Urba\'nski \cite{PrzUrb}
proved the inequality~(\ref{dimension}) for the case
of the Erd\"os measure $\mu$ and Pisot number $\lambda$
(and the equality in (\ref{dimension}) was shown by S.~Lalley \cite{Lalley}
-- see Example~\ref{Erd} below).

V.~Kaimanovich \cite{Kai2}
obtained a similar result for the Hausdorff dimension
of the harmonic measure on trees with applications to certain classes
of random walks.
R.~Lyons \cite{Lyons} has a result analogous to Corollary~\ref{cor-1}
in the context of random continued fractions. S.~Pincus \cite{Pincus}
has related results in the context of mappings on the line and
$2\times2$ matrices in the plane. K.~Simon, B.~Solomyak and
M.~Urba\'nski \cite{SSU} have a theorem similar to Corollary~\ref{cor-1}
in the context of parabolic iterated function systems
on the real line. Moreover, they were able to establish
certain parameter values of their system for which the measure $\mu$
is absolutely continuous a.e.
\end{rmk}

\begin{exam} Let us give a simple example.
Suppose $f_1(x)=2x+1,f_2(x)=\frac1{16} x+1$ chosen
with probabilities $p_1=p_2=\frac12$. Then
$\chi(\Phi)=-\frac32\log 2< h(\ov p)=-\log 2$ and hence by
Corollary~\ref{cor-1}, the invariant measure $\mu$ is singular with respect
to Lebesgue measure, and $\dim_H(\mu)\le\frac23$.
However it is easy to show that the support of
the invariant measure is the interval $[1,\infty)$.
Note that a more detailed analysis shows that since
$f_1f_2f_1^3f_2f_1=f_2f_1^5f_2$, we have $\th<1.9836$, whence
by the formula~(\ref{dim-th}), $\dim_H(\mu)\le\frac23\log_2\th<0.6588$.
For more examples see Section~\ref{examples}.
\end{exam}

\noindent {\bf Proof of Theorem~\ref{thm-2}:}

We let $B(x,r)$ denote the ball of radius $r$ about the
point $x\in \R^d$. It is known (see \cite[Page 171]{Falconer})  that
for any Borel measure $\mu$ on $\R^d$,
\begin{equation}
\dim_H(\mu)=\mu\mbox{-ess}
\sup \left\{\liminf_{r\rightarrow 0} \frac{\log \mu B(x,r)}{\log r}\right\}.
\label{ess-sup}
\end{equation}
Let $B_{\om}:=B(\phi(\om),1)$ and
$f_{\om}^{(n)}:=f_{\om_0}\dots f_{\om_{n-1}}$. Our goal is to show that
\begin{equation}
\mu f^{(n)}_{\om^*}(B_{\om^*})\ge
\mu(B_{\om^*})\cdot\nu\mathcal{E}_n(\om^*).
\label{mainest}
\end{equation}
We have
$\mu f^{(n)}_{\om^*}(B_{\om^*})=\nu(\phi^{-1}f^{(n)}_{\om^*}(B_{\om^*}))$
and
\[
\aligned
\phi^{-1}f^{(n)}_{\om^*}(B^*)=\{&\om:\lim_{k\to\infty}
f_{\om_0}\dots f_{\om_n}\dots f_{\om_{n+k-1}}(x_0)\in
f^{(n)}_{\om^*}(B^*)\} \\
\supset \{&\om:(\om_0,\dots,\om_{n-1})\sim(\om^*_0,\dots,\om^*_{n-1}),\\
&\lim_{k\to\infty} f_{\om_n}\dots f_{\om_{n+k-1}}(x_0)\in B_{\om^*} \},
\endaligned
\]
whence by the fact that $\nu$ is a product measure
(\ref{mainest}) follows. Hence by Lemma~\ref{shannon}
for any fixed $\delta>0$ for
$\nu$-a.e. $\om^*$ for all sufficiently large $n$,
$$
\mu f^{(n)}_{\om^*} (B_{\om^*})\ge\mu(B_{\om^*})\th^{-n(H_\mu+\delta)}.
$$
We define $\gamma=\gamma(\Phi):=\exp \chi(\Phi)$,
fix $\delta>0$ sufficiently small that
$0<\gamma-\delta<\gamma + \delta <1$ and
define the sets $G^1_N$ and $G^2_N$ as follows:
\[
\begin{aligned}
G^1_N &:= \{ \om\in\Om: \forall n> N,
(\gamma-\delta)^n < \|f_{\om_0}\dots f_{\om_{n-1}}\|<(\gamma+\delta)^n\},
\\
G^2_N&:=\{\omega\in\Om :
\forall n> N,
\mu f^{(n)}_{\om^*} (B^*)\ge\mu(B_{\om^*})\th^{-n(H_\mu+\delta)}\}.
\end{aligned}
\]
Let $G_N=G^1_N \cap G^2_N$. We may choose $N$ sufficiently large
that $\nu (G_N)>\frac34$.

Note that $\nu$-a.e. $\omega\in \Omega$ has the
property that $\phi(\omega) \in \mbox{ supp}(\mu)$.
Let $\om$ be such a sequence and define
$A(N):=\{\om : \mu(B_\om)>\frac1N\}$. Since $\mu(B_\om)>0$,
we have $\nu\left(\bigcup_{N=1}^\infty A(N)\right)=1$. As
$A(N+1)\subset A(N)$, we have $\lim_N \nu(A(N))=1$ as well.

Hence we may fix $\alpha>0$ sufficiently small such
that $\nu\{\om : \mu(B_\omega)>\alpha \}>\frac34$.
Define
$$
B:=\{ \omega \in G_N \mid \mu(B_\omega)>\alpha \}.
$$
Then $\nu (B)>\frac12$ and by the fact that $\sigma$ preserves the measure
$\nu$, we have $\nu(\sigma^n (B))>\frac12$ for any $n>0$.

We claim that for any $n>N$ and any $x\in \phi (\sigma^n (B))$,
\begin{equation}
\frac{\mu B(x,r)}{\L_d B(x,r)}\ge C'(d)\alpha
\th^{-n(H_\mu+\delta)}(\gamma+\delta)^{-dn}\label{muoverm}
\end{equation}
for some $r>0$ (here $C'(d)$ is a constant which depends upon the
dimension $d$ and $\L_d$ is $d$-dimensional Lebesgue measure).
To prove (\ref{muoverm}), for $n>N$ and $\om^* \in B$ we let
$x=\phi(\sigma^n \omega^*)$  and $r=\|f^{(n)}_{\om^*}\|$.
Since $f^{(n)}_{\om^*}(B^*)\subset B(x,r)$, we have
$$
\mu B(x,r)\ge\alpha\th^{-n(H_\mu+\delta)}.
$$
To estimate $\L_d B(x,r)$, we note that since $\om^*\in G_N$, we have
$r<(\gamma+\delta)^n$, whence
$$
\L_d B(x,r)\le C_d(\gamma+\delta)^{dn},
$$
where $C_d$ is the volume of the unit ball in $\R^d$.
This proves (\ref{muoverm}) with $C'(d)=1/C_d$.

Since $(\Omega,\sigma,\nu)$  is ergodic,
for $\nu$-a.e. $\omega$ we have $\omega \in \sigma^n (B)$
for infinitely many integers $n$. Hence for $\mu$-a.e. $x\in\mbox{supp}(\mu)$
we have $x\in \phi(\sigma^n B)$ infinitely often.

This establishes the fact that for a $\mu$-generic
$x\in\mbox{supp}(\mu)$ there exists a subsequence $r_n \rightarrow 0$
such that,
\[
\frac{\mu B(x,r_n)}{\L_d B(x,r_n)}\ge C'(d)\alpha
\th^{-n(H_\mu+\delta)}(\gamma+\delta)^{-dn},
\]
which is equivalent to
\[
\frac{\mu B(x,r_n)}{r_n^d}\ge\alpha
\th^{-n(H_\mu+\delta)}(\gamma+\delta)^{-dn},
\]
since $C'(d)=1/C_d$ and $\L_d B(x,r_n)=C_d r_n^d$. Taking logarithms
and dividing by $\log r_n$, we have
\[
\frac{\log\mu B(x,r_n)}{\log r_n} -d\le
\frac{\log\alpha}{\log r_n}-
\frac{n}{\log r_n}\left((H_\mu+\delta)\log\th+d\log(\gamma+\delta)\right).
\]
Since $x=\phi(\sigma^n\om^*)$, where $\om^*\in G_N$, we have
$(\gamma-\delta)^n\le r_n \le (\gamma+\delta)^n$, whence it follows that
for $\mu$-a.e. $x$,
\[
\begin{aligned}
\liminf_{r \rightarrow 0}\frac{\log \mu B(x,r)}{\log r}
&\le
\liminf_{r_n \rightarrow 0}\frac{\log \mu B(x,r_n)}{\log r_n}\\
&\le
d-\frac{(H_\mu+\delta)\log\th}{\log\gamma}-
d\ \frac{\log(\gamma+\delta)}{\log\gamma}.
\end{aligned}
\]
Since $\delta>0$ may be taken arbitrarily small and $\log\gamma=\chi(\Phi)$,
we finally obtain
$$
\liminf_{r\to0}\frac{\log\mu B(x,r)}{\log r}\le
-\frac{H_\mu\log\th}{\chi(\Phi)},
$$
and by (\ref{ess-sup}) inequality~(\ref{dimension}) holds,
which completes the proof.\qed

\section{Examples}\label{examples}

We are going to consider several examples, all of which
are affine IFS.

\begin{exam}\label{Erd} (Bernoulli convolutions).
Put $\Om:=\prod_0^\infty\{0,1\}$
and let $\la>1,\ f_0(x)=\la^{-1}x,f_1(x)=\la^{-1}x+1-\la^{-1},
p_1=p_2=\frac12$ (see Introduction). In this case $\chi(\Phi)=-\log\la$, and
$$
f_{\om_0}\circ\dots\circ f_{\om_{n-1}}(0)=
(1-\la^{-1})\sum_{k=0}^{n-1}\om_k\la^{-k},
$$
thus,
$$
\phi(\om)=(1-\la^{-1})\sum_{k=0}^\infty\om_k\la^{-k}.
$$
Hence $(\om_0,\dots,\om_{n-1})\sim(\om'_0,\dots,\om'_{n-1})$ iff
$\sum_0^{n-1}\om_k\la^{-k}=\sum_0^{n-1}\om'_k\la^{-k}$. Assume first
$\la>2$. Then $\mbox{supp}(\mu)$ is known to be a Cantor set, the semigroup
$G^+$ is obviously free (as $f_0([0,1])\cap f_1([0,1])=\emptyset$),
and from Corollary~\ref{cor-1} it follows
\begin{equation}
\dim_H(\mu)\le\frac{\log 2}{\log\la}<1.
\label{log2/logla}
\end{equation}

If $1<\la<2$, then $\mbox{supp}(\mu)=[0,1]$;
if $\la$ is transcendental, it is easy to see that $G^+=\F_2^+$, whence
$\th=2$ and $H_\mu=1$. Hence Corollary~\ref{cor-1} again gives us the
estimate~(\ref{log2/logla}),
which is unfortunately useless, as $\la<2$. However, in some cases of
algebraic $\la$ Theorem~\ref{thm-2} can be used more efficiently.

More specifically, assume $\la$ to be a {\em Pisot number}, i.e., an
algebraic integer greater than 1 whose conjugates have moduli less than
1. The famous Separation Lemma due to A.~Garsia
\cite{Garsia} states that there exists a constant $C=C(\la)>0$ such that
if $\sum_0^n\om_k\la^{-k}\neq\sum_0^n\om'_k\la^{-k}$, then
$\left|\sum_0^n(\om_k-\om_k')\la^{-k}\right|\ge C\la^{-n}$. Hence
it is easy to see that $\th=\la$, and from (\ref{dimension}) follows
$\dim_H(\mu)\le H_\mu$.

In work by S.~Lalley \cite{Lalley} the Separation Lemma was used to show
that in fact
$$
\dim_H(\mu)=H_\mu<1.
$$
The most transparent subcase is $\la=\frac{1+\sqrt5}2$. It was studied
in several papers (see references in \cite{SidVer});
in particular, for this $\la$ we have
$G^+=\langle a,b\mid ab^2=ba^2\rangle$ and $\dim_H(\mu)=H_\mu=0.995713\dots$
(this numerical result is due to J.~C.~Alexander and D.~Zagier
\cite{AlZa}). Besides, the measure $\mu$ was shown to be
quasi-invariant under the {\em $\beta$-shift} (for $\beta=\la$)
$\tau_\la:[0,1)\to[0,1)$
defined by the formula $\tau_\la x=\{\la x\}$ and the corresponding density
is also known (see \cite{SidVer}).
\end{exam}
Similar results hold for a more general Bernoulli convolution
$\mu=B_\la(p,1-p)$, i.e., the one for which the probability of taking
$f_0$ is $p\in(0,1)$.

We believe the techniques of \cite[Proposition~4]{Lalley}
can be used to show that the equality holds in a more general situation.
Let us formulate the corresponding conjecture; put
as above, $x_n(\om):=f_{\om_0}\dots f_{\om_{n-1}}(x_0)$ and
$\gamma=\exp\chi(\Phi)$.

\smallskip\noindent
{\bf Conjecture.} Suppose we have an affine IFS on the real line (i.e.,
$f_i(x)=\la_ix+b_i$) and $|\la_i|\ge1$ for at least one $i\in\{1,\dots,m\}$;
assume the following Weak Separation Condition to be satisfied (we borrow
this term from \cite{LauNgaiRao}):
for $\nu^2$-a.e. pair $(\om,\om')\in\Omega^2$ and arbitrary $\delta>0$,
\begin{equation}
|x_n(\om)-x_n(\om')|\ge\mbox{const}\cdot(\gamma-\delta)^n
\label{WSC}
\end{equation}
whenever $x_n(\om)\neq x_n(\om')$.
We {\em conjecture} that the inequality~(\ref{dimension}) in this context
is actually an equality. We state without proof that ({\ref{WSC}) does hold
in the framework of Example~\ref{simplest} with $\la$ being a Pisot
number (see below).

Suppose $\Phi$ is an affine IFS in $\R$.
If $\Phi$ is not uniformly contracting (i.e., there exists
$i$ such that $|\la_i|\ge1$), then
by the result from \cite{mBb88} mentioned above, the support of
the invariant measure in this case is either a single point or
$\R$ or $[d,+\infty)$ or finally $(-\infty,d]$ for some $d\in\R$.
We may rule out the first case. The fact that $\mbox{supp}(\mu)$ is
connected makes the problem about the fine structure of $\mu$ nontrivial.

\begin{exam}\label{simplest}
Let $\la>1$ and
\begin{equation}
f_1(x)=\la^{-1}x,\ f_2(x)=x+1 \textrm{ and }
p_1=p_2=\frac12.
\label{basic}
\end{equation}
The support of $\mu=\mu(\la)$
is $[0,+\infty)$, and $\chi(\Phi)=-\frac12\log\la$. Hence by
Corollary~\ref{cor-1}, for $\la>4$ the measure $\mu$ is singular,
and $\dim_H(\mu)\le \frac{2\log 2}{\log\la}<1$.

We claim that for any transcendental $\la$ the
semigroup $G^+$ is free. A trivial induction argument shows that
$$
f_1^{n_1}f_2^{k_1}\dots f_1^{n_s}f_2^{k_s}(x)=\la^{-\sum_1^s n_j}x+
\sum_{j=1}^s k_j\la^{-\sum_{i=1}^j n_i},
$$
whence if $\la$ is not algebraic,
$$
f_1^{n_1}f_2^{k_1}\dots f_1^{n_s}f_2^{k_s}=
f_1^{n'_1}f_2^{k'_1}\dots f_1^{n'_s}f_2^{k'_s}
$$
implies $n_j\equiv n'_j,\ k_j\equiv k'_j,\ j=1,\dots,s$. Hence
$\#D_n=2^n,\ \th=2$ and $G^+=\F_2^+$.

It is worth noting that for certain particular values of $\la\in(1,4)$
the measure $\mu$ is nonetheless singular (similarly to the
Bernoulli convolutions).
Namely, since $\mu$ is invariant under the IFS~$\Phi$,
we have the following self-similar relation for its Fourier transform:
$$
\wh\mu(x)=\frac12\ \wh\mu(\la^{-1}x)+\frac12\ e^{ix}\wh\mu(x),
$$
whence
\begin{equation}
\wh\mu(x)=\prod_{n=0}^\infty\frac1{2-\exp(i\la^{-n}x)}.
\label{fourier}
\end{equation}
Assume again $\la$ to be a Pisot number. Then as is well known, the distance
to the nearest integer for $\la^n$ tends to 0 at exponential rate. Following
the line of the proof of the classical work \cite{Erdos} (see also
\cite{mBb88} for the case $\la=2$), we can consider
the subsequence $x_n=2\pi\la^n$ and show that $\wh\mu(x_n)\not\to0$ as
$n\to+\infty$, which implies that the Riemann-Lebesgue Lemma is not satisfied,
whence $\mu$ cannot be absolutely continuous. Therefore, it is singular
by the Law of Pure Types.

Thus, we proved that $\mu$ is singular for $\la\ge4$ (as 4 is a Pisot
number); at the same time it is singular for
an infinite number of parameters $\la\in(1,4)$ as well.
It is an open question whether its Hausdorff
dimension is less than 1 for a Pisot number $\la$ (for the Bernoulli
convolutions it is true, see above).

It is worth mentioning that in this example the stationary measure has
an ``arithmetic" interpretation as well. Namely, let
$\Sigma=\prod_1^\infty\Z_+$ and $\xi$ denote the stationary product measure
on $\Sigma$ with the following geometric
distribution: $\xi(\eps_n=k)=2^{-k-1},\
k=0,1,\dots$ Let $L_\la:\Sigma\to\R_+$ be defined as follows:
$$
L_\la(\eps_1,\eps_2,\dots)=\sum_{n=1}^\infty \eps_n\la^{-n}
$$
(it is obvious that $L_\la$ is well defined for $\xi$-a.e. $\eps\in\Sigma$).
Then from (\ref{fourier}) it follows
$$
\mu=\xi\circ L_\la^{-1},
$$
as $\wh\xi(x)=\frac12+\frac14e^{ix}+\frac18e^{2ix}+\dots=\frac1{2-e^{ix}}$.
Thus, the essential difference with the case of Bernoulli convolutions is
that the set of ``digits" here is infinite.

When this paper was in preparation, Y.~Peres suggested the following claim.

\begin{lemma}
For $\L_1$-a.e. $\la\in(1,3\cdot 2^{-2/3})$ the measure $\mu$ is absolutely
continuous.
\end{lemma}
\noindent {\bf Proof.} Since
$$
\frac1{2-e^{ix}}=\frac{2+e^{ix}}{4-e^{2ix}}=
\frac{\frac23+\frac13e^{ix}}{\frac43-\frac13e^{2ix}},
$$
by (\ref{fourier}) the measure $\mu$ is the convolution of the Bernoulli
measure $B_\la(2/3,1/3)$ and the stationary measure for the IFS
$$
\begin{aligned}
g_1(x)&=\la^{-2}x,\\
g_2(x)&=x+1
\end{aligned}
$$
with $p_1=\frac14,\ p_2=\frac34$. In \cite{PeresSol2} it was shown that
for any $p\in[1/3,2/3]$ the Bernoulli measure is absolutely continuous
for a.e. $\la\in(1,p^{-p}(1-p)^{-1+p})$. Hence
the measure $B_\la(2/3,1/3)$ will be absolutely continuous for a.e.
$\la\in(1,3\cdot 2^{-2/3})$ and so will be $\mu$.\qed

\begin{rmk} Since
\[
\frac1{1-z}=\prod_{n=0}^\infty \left(1+z^{2^n}\right),\quad |z|<1,
\]
it is easy to deduce that
$$
\frac1{2-e^{ix}}=\prod_{n=0}^\infty\left(\frac{2^{2^n}}{2^{2^n}+1}+
\frac1{2^{2^n}+1}\,\exp\left(2^nix\right)\right),
$$
and by (\ref{fourier})
\begin{equation}
\mu=B_\la\left(\frac23,\,\frac13\right)*
B_{\la,2}\left(\frac45,\,\frac15\right)*
\dots*B_{\la,{2^n}}\left(\frac{2^{2^n}}{2^{2^n}+1},\frac1{2^{2^n}+1}\right)*
\dots,
\label{conv}
\end{equation}
where $B_{\la,t}(p,1-p)$ is the infinite convolution of two-point
discrete measures whose ``basic" distribution is supported on the points
0 and $t$ with probabilities $p$ and $1-p$ respectively (hence
$B_\la=B_{\la,1}$).
\end{rmk}
Summing up, we have the following

\begin{prop} For the IFS~(\ref{basic}) the following properties
are satisfied:
\begin{enumerate}
\item For $\la\ge4$ the measure $\mu$ is singular, and
$\dim_H(\mu)\le\frac{2\log2}{\log\la}$;
\item for $\L_1$-a.e. $\la\in(1,3\cdot 2^{-2/3})$ it is absolutely
continuous;
\item for any Pisot $\la$ it is singular.
\end{enumerate}
\end{prop}

A natural question to ask is what happens between
$3\cdot 2^{-2/3}\approx 1.8899$ and 4. Note that if $\la>3\cdot 2^{-2/3}$,
then {\em all} convolutions in (\ref{conv}) are singular.
Nonetheless, we {\em conjecture} that for $\L_1$-a.e.
$\la\in(3\cdot 2^{-2/3},4)$ the measure $\mu$ will be absolutely continuous
as well.
\end{exam}

\begin{exam} Let $\la>1,\ f_1(x)=\la x+1,\ f_2(x)=\la^{-2}x$,
and $p_1=p_2=\frac12$. Here $\chi(\Phi)=-\frac12\log\la$.
Let for the simplicity of notation $a=f_1,\ b=f_2$. Since
\begin{equation}
ababa=ba^3b,\label{ababa}
\end{equation}
we have $d_{n+1}\le2d_n-d_{n-4}$, whence
$\theta<1.9277$. Thus, from the estimate~(\ref{dim-th}) it follows
that the measure $\mu$ is
singular at least for $\la>1.9277^2\approx 3.716$. However, it is
clear that the actual estimate must be even sharper, because
in fact there are infinitely many relations between $a$ and $b$.

Namely, from general considerations it follows that $G^+$
will be one and the same for any transcendental $\la$. We do not know whether
$G^+$ is finitely presented but at least  written in the generators
$a,b$ it is not. For example, for any $n=2,3,\dots$ and any
$k=1,2,\dots,n-1$ we have in addition to the relation~(\ref{ababa}),
$$
a^{2k}b^n a^{2(n-k)}=b^{n-k}a^{2n}b^k
$$
(these are a direct consequence of the fact that $a^{2k}b^k$ and $a^{2n}b^n$
always commute). These relations are independent, i.e., none of them
is a consequence of any other one. At the same time, there are relations
that cannot be deduced from any of those described above; for instance,
$abab^2a^3=b^2a^5b$. There are indications that actually $\th$ is at least
less than 1.7.
\end{exam}

As far as we are concerned, there are no general results on the structure
of $\mbox{supp}(\mu)$ in the case of higher dimensions. We are going
to present an example of a family of IFS for $d=2$ such that
$\mbox{supp}(\mu)=\R^2$, whereas $\mu$ is singular; at the same time
this system does not ``split" into one-dimensional actions. In a
certain sense the following example is a two-dimensional generalization
of Example~\ref{simplest}.

\begin{exam} Let $\al$ be a real number such that
$\al/\pi$ is irrational and let
$R_\al$ denote the rotation of $\R^2$ by the angle $\al$, i.e.,
$$
R_\al=\left(\begin{matrix} \cos\al&-\sin\al \\ \sin\al&\cos\al
\end{matrix}\right).
$$
Let $\la>1$ and the one-parameter family of IFS $\Phi_\la$ be
defined as follows:
$$
\begin{aligned}
f_1(x)&=A_\la^{-1}x,\\
f_2(x)&=x+\ov e,
\end{aligned}
$$
where $A_\la=\la R_\al$ and $\ov e\in S^1$ is fixed. As above,
we assume $p_1=p_2=\frac12$.

\begin{prop}
\begin{enumerate}
\item For any $\la>1$ the IFS $\Phi_\la$ contracts on average and the support
of the invariant measure $\mu_\la$ is full, i.e.,
\begin{equation}
\mbox{supp}(\mu_\la)=\R^2.
\label{full}
\end{equation}
\item For $\la>2$ the measure $\mu_\la$ is singular, and
\begin{equation}
\dim_H(\mu_\la)<\frac{2\log 2}{\log\la}<2.
\label{dim<d}
\end{equation}
\end{enumerate}
\end{prop}
\noindent {\bf Proof:} By the same reason as in the previous examples,
we have $\chi(\Phi_\la)=-\frac12\log\la$, whence for any $\la>1$ the
system contracts on average. From Corollary~\ref{cor-1} it follows that
$\la>2$ implies the singularity of $\mu_\la$ together with (\ref{dim<d}).

The most delicate part of the proposition is the relation~(\ref{full}).
Let us prove it. Assume $\M:=\mbox{supp}(\mu_\la)\neq\R^2$; then there
exists a disc $B(x,\delta)$ whose intersection with $\M$ is empty.
Hence by definition, $f_1^{-n}B(x,\delta)\cap\M=\emptyset$ for any
$n\ge1$. We have
$$
f_1^{-n}B(x,\delta)=B(y_n,\la^n\delta),
$$
where $y_n=f_1^{-n}(x)=\la^n R_\al^n(x)$. Since $\al/2\pi$ is irrational,
the rotation $R_\al:S^1\to S^1$
is {\em minimal}, i.e., the orbit of every point is dense in $S^1$
(see, e.g., \cite{CorFoSin}). We apply this claim to the circle
of radius
$\|x\|$. Thus, for any $\eps>0$ there exists
$n\in\mathbb{N}$ such that
\begin{equation}
\|R_\al^n(x)-\|x\|\ov e\ \|<\eps.
\label{<eps}
\end{equation}
Fix $r>1,\ \eps=\delta/2$ and $n$ large enough to satisfy
$\la^n\ge 2r/\delta$ together with (\ref{<eps}). Let $z=\la^n\|x\|\ov e$;
we claim that
$$
B(z,r)\subset B(y_n,\la^n\delta).
$$
Indeed, let $y\in B(z,r)$, i.e., $\|y-z\|\le r$.
Hence
$$
\begin{aligned}
\|y-y_n\|&\le \|y-z\|+\|z-y_n\| \\
         &\le r+\la^n\eps=r+\frac12\la^n\delta<\la^n\delta.
\end{aligned}
$$
Hence $B(z,r)\cap\M=\emptyset$, and $f_2^{-k}B(z,r)\cap\M=\emptyset$
as well for any $k\ge0$. Since $\|\ov e\|=1$ and $z$ belongs to
the half-line $\{t\ov e, t\ge0\}$, there exists $k\ge0$ such that
$f_2^{-k}B(z,r)\supset B(0,r-1)$. Hence for any
$r>1,\ B(0,r-1)\cap\M=\emptyset$, which means $\M=\emptyset$.
The proposition is proven.\qed

\end{exam}


\begin{thebibliography}{999}

\bibitem{AlZa}
\newblock Alexander J. C. and Zagier D. (1991).
\newblock {\em The entropy of a certain infinitely convolved Bernoulli
measure},
\newblock J. London Math. Soc. {\bf 44}, 121--134.

\bibitem{Arnold&Crauel}
\newblock Arnold L. and Crauel H. (1990).
\newblock  {\em Iterated Function Systems and Multiplicative Ergodic Theory},
\newblock Diffusion processes and related problems in analysis, Vol. II
(Charlotte, NC), 283--305.

\bibitem{mBb88}
\newblock Barnsley M. F. and Elton J. H. (1988).
\newblock {\em A new class of Markov processes for image encoding},
\newblock Adv. Appl. Prob., {\bf 20}, 14--32.

\bibitem{BHN}
\newblock Broomhead D., Hadjiloucas D. and Nicol M. (1999).
\newblock {\em Random and deterministic perturbation
of a class of skew-product systems},
\newblock Dynamics and Stability of Systems, {\bf 14}, No.~2, 115--128.

\bibitem{Campbell}
\newblock Campbell K. M. (1997).
\newblock {\em Observational noise in skew product systems},
\newblock {Physica D} {\bf 107}, 43--56.

\bibitem{CorFoSin}
\newblock Cornfeld I. P., Fomin S. V. and Sinai Ya. G. (1982).
\newblock {\bf Ergodic theory},
\newblock New York-Heidelberg-Berlin: Springer-Verlag.

\bibitem{Der}
\newblock Derriennic Y. (1980).
\newblock {\em Quelques applications du th\'eor\`eme ergodique
sous-additif},
\newblock Ast\'erisque {\bf 74}, 183--201.

\bibitem{DiFr}
\newblock Diaconis P. and Freedman D. (1999).
\newblock {\em Iterated random functions},
\newblock SIAM Review {\bf 41}, 45--76.

\bibitem{Dubins-Freedman}
\newblock Dubins L. E. and Freedman D. A. (1966)
\newblock {\em Invariant probabilities for certain Markov processes},
\newblock Ann. Math. Stat. {\bf 32}, 837--848.


\bibitem{El2} Elton J. H. (1990).
\newblock {\em A multiplicative ergodic theorem for Lipschitz maps},
\newblock {Stochastic Processes and their Applications} {\bf 34},
39--47.

\bibitem{Erdos}
\newblock Erd\"os P. (1939).
\newblock {\em On a family of symmetric Bernoulli convolutions},
\newblock Amer. J. Math. {\bf 61}, 974--975.

\bibitem{Falconer} Falconer K. (1997).
\newblock {\bf Techniques in Fractal Geometry},
\newblock John Wiley and Sons.

\bibitem{Garsia}
\newblock Garsia A. (1962).
\newblock {\em Arithmetic properties of Bernoulli convolutions},
\newblock Trans. Amer. Math. Soc. {\bf 102}, 409--432.

\bibitem{Kai}
\newblock Kaimanovich V. (1997)
\newblock The Poisson formula for groups with hyperbolic properties,
\newblock preprint,
\newblock WWW: http://www.maths.univ-rennes1.fr/\symbol{126}kaimanov

\bibitem{Kai2}
\newblock Kaimanovich V. (1998)
\newblock {\em Hausdorff dimension of the harmonic measure on trees},
\newblock Ergodic Theory Dynam. Systems {\bf 18}, 631--660.

\bibitem{KaiVer}
\newblock Kaimanovich V. and Vershik A. (1983).
\newblock {\em Random walks on discrete group: boundary and entropy},
\newblock {Ann. Prob.} {\bf 11}, 457--490.

\bibitem{Lalley}
\newblock Lalley S. (1998).
\newblock{\em Random series in powers of algebraic numbers: Hausdorff
dimension of the limit distribution},
\newblock J. London Math. Soc. {\bf 57}, 629--654.

\bibitem{LauNgaiRao}
\newblock Lau K.-S., Ngai S.-M. and Rao H. (1999)
\newblock {\em Iterated function systems with overlaps and self-similar
measures}, to appear in J. London Math. Soc.

\bibitem{Lyons}
\newblock Lyons R. (1999).
\newblock {\em Singularity of some random continued fractions},
\newblock  J. Theor. Prob., to appear.

\bibitem{PeScSo}
\newblock Peres Y., Schlag W. and Solomyak B. (2000).
\newblock {\em Sixty years of Bernoulli convolutions},
Fractals and Stochastics II, (C. Bandt, S. Graf and M. Zaehle,
eds.), Progress in Probability {\bf 46}, 39--65, Birkhauser.

\bibitem{PeresSol2}
\newblock Peres Y. and Solomyak B. (1998)
\newblock {\em Self-similar measures and intersections of Cantor sets},
\newblock Trans. Amer. Math. Soc. {\bf 350}, 4065--4087.

\bibitem{Pincus}
\newblock Pincus S. (1994).
\newblock  {\em Singular stationary measures are not always fractal},
\newblock J. Theor. Prob., {\bf 7}, 199--208.


\bibitem{PrzUrb}
\newblock Przytycki F. and Urba\'nski M. (1989)
\newblock {\em On the Hausdorff dimension of some fractal sets},
\newblock Studia Math. {\bf 93}, 155--186.

\bibitem{SidVer}
\newblock Sidorov N. and Vershik A. (1998).
\newblock {\em Ergodic properties of Erd\"os measure, the entropy of
the goldenshift, and related problems},
\newblock Monatsh. Math. {\bf 126}, 215--261.

\bibitem{SSU} Simon K., Solomyak B. and Urba\'nski M. (1999).
\newblock {\em Invariant measures
for parabolic IFS with overlaps and random continued fractions}, preprint

\bibitem{Solomyak} Solomyak B. (1995).
\newblock {\em On the random series
$\sum\pm \lambda^i$ (an Erd\"os problem)},
\newblock Annals of Math. {\bf 142}, 611--625.

\bibitem{Stark1}
\newblock Stark J. (1999).
\newblock {\em Regularity of invariant graphs for forced systems},
\newblock {Ergodic Theory  Dynam. Systems}, {\bf 19}, 155--199.

\end{thebibliography}
\end{document}